\documentclass[12pt]{article}
\usepackage{amsmath,amssymb,amsbsy}
\usepackage{graphicx}
\newtheorem{thm}{Theorem}

\newtheorem{cor}[thm]{Corollary}

\newcommand{\prf}{\noindent{\it Proof.}\ }
\newcommand{\qed}{\hfill \rule{.1in}{.1in}}
\def\mod{\mathop{\rm mod}\nolimits}

\newcommand{\ent}{\mathbb N}

\begin{document}
\title{Super edge-graceful paths}
\author{Sylwia Cichacz\thanks{The work was supported by Fulbright Scholarship nr 15072441.} \thanks{cichacz@agh.edu.pl} \\
{\small AGH University of Science and Technology}\\
{\small and}\\
{\small University of Minnesota Duluth}\\
Dalibor Froncek\\
{\small University of Minnesota Duluth}\\
Wenjie Xu\\
{\small University of Minnesota Duluth}}

\maketitle

\date

\begin{abstract}
A graph $G(V,E)$ of order $|V|=p$ and size $|E|=q$ is called super
edge-graceful if there is a bijection $f$ from $E$ to $\{0,\pm 1,\pm
2,\ldots,\pm \frac{q-1}{2}\}$ when $q$ is odd and from $E$ to $\{\pm
1,\pm 2,\ldots,\pm \frac{q}{2}\}$ when $q$ is even such that the
induced vertex labeling $f^*$ defined by $f^*(x) = \sum_{xy\in
E(G)}f(xy)$ over all edges $xy$ is a bijection from $V$ to $\{0,\pm
1,\pm 2\ldots,\pm \frac{p-1}{2}\}$ when $p$ is odd and from $V$ to
$\{\pm 1,\pm 2,\ldots,\pm \frac{p}{2}\}$ when $p$ is even. \\
\indent We prove that all paths $P_n$ except $P_2$ and $P_4$ are
super edge-graceful.
\end{abstract}
\section{Introduction}

All graphs considered in this paper are simple, finite and
undirected. We use  standard terminology and notation of graph
theory.\\
\indent Let $G=(V,E)$ be a graph with $p$ vertices and $q$ edges. A
vertex labeling of a graph $G$ is a function from $V(G)$ into
$\ent$. A. Rosa \cite{R} introduced the graceful graph labeling. A
graph $G$ is \emph{graceful} if there exists  an injection from the
vertices of $G$ to the set $\{0,\ldots,q\}$ such that, when each
edge $xy$ is assigned the label $|f(x)-f(y)|$, the resulting
edge labels are distinct.\\
\indent The edge-graceful labeling was introduced by S.P. Lo
\cite{L}. A graph $G$ is \emph{edge-graceful} if the edges can be
labeled by $1,2,\ldots,q$ such that the vertex sums are distinct
($\mod p$). A necessary condition for a graph with $p$ vertices and
$q$ edges to be edge-graceful is that $q(q +1) \equiv
\frac{p(p-1)}{2} (\mod p)$.

J. Mitchem and A. Simoson \cite{MS} defined super edge-graceful
labeling which is a stronger concept than edge-graceful for some
classes of graphs.\\
Define an edge labeling as a bijection $$f : E(G) \rightarrow
\{0,\pm 1,\pm 2,\ldots,\pm \frac{q - 1}{2}\}\,\, \mathrm{ for }\,\,
q \, \, \mathrm{ odd}$$ or
$$f : E(G) \rightarrow \{\pm 1,\pm 2,\ldots,\pm \frac{q}{2}\} \,\, \mathrm{for}\,\, q \,\, \mathrm{even}.$$
For every vertex $x \in V (G)$, define the induced vertex labeling
of $x$ as $f^*(x) = \sum_{ xy\in E(G)}f(xy)$. If $f^*$ is a
bijection
$$f^* : V(G) \rightarrow \{0,\pm 1,\pm 2,\ldots,\pm \frac{p -
1}{2}\}\textrm{ for } p \textrm{ odd }$$ or $$f^* : V(G) \rightarrow
\{\pm 1,\pm 2,\ldots,\pm p\} \textrm{ for } p \textrm{ even},$$ then
the labeling $f$ is
\emph{super edge-graceful}. \\

S.-M. Lee and Y.-S. Ho showed that all trees of odd order with three
even vertices are super edge-graceful \cite{LH}. In \cite{CLGS}
P.-T. Chung, S.-M. Lee, W.-Y. Gao, and K. Schaffer asked which paths
are super edge-graceful. In this paper we show that all paths $P_n$
except $P_2$ and $P_4$ are super edge-graceful.

\section{Super edge-gracefulness of $P_n$}
\begin{figure}[h]
\begin{center}
\includegraphics[width=11.5cm]{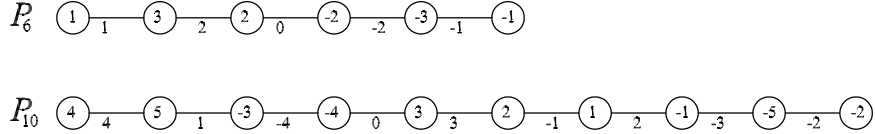}
\caption{A super edge-graceful labeling $P_{6}$ and
$P_{10}$.}\label{6and10}
\end{center}
\end{figure}

Let $P_m=x_1,x_2,\ldots,x_m$ be a path with an edge labeling $f$ and
$P_m'=x_1',x_2',\ldots,x_m'$ be a path with an edge labeling $f'$.
If for every $i$, $1 \leqslant i \leqslant m-1$ we have
$f(x_ix_{i+1})=-f'(x_ix_{i+1})$ then the labeling $f$ is called
\emph{inverse} of $f'$.
\begin{thm} The path $P_n$ is super edge-graceful unless $n
= 2,4$.\label{CFX}
\end{thm}
\prf It is obvious that $P_2$ is not super edge-graceful. $P_4$ is
not super edge-graceful since the edge label set is $\{0, -1, 1\}$
and the vertex set is $\{-2, -1, 1, 2\}$, but no two edge labels
will sum up to $2$ or $-2$. A labeling of $P_3$ is trivial. We label
the edges along $P_6$ by $(1,2,0,-2,-1)$, whereas along $P_{10}$
by $(4,1,-4,0,3,-1,2,-3,-2)$ (see Figure~\ref{6and10}).\\

Assume from now on that $n\geqslant 5$ and $n \neq 6,10$. The basic
idea of our proof is to consider a path $P_n$ as a union of paths
$Q_{1}$ and $Q_{2}$ joined by an edge with
label $0$. We consider the following cases (for the sake of completeness we will include cases for odd paths $P_n$ for $n\geqslant 5$, which were proved in \cite{LH}):\\

{\noindent \emph{Case 1.}} $n\equiv 1 (\mod 4)$.

It follows that $n=4k+1$ where $k$ is a positive integer. Then we
can find a super edge-graceful labeling as follows:\\
\begin{figure}[h]
\begin{center}
\includegraphics[width=13.5cm]{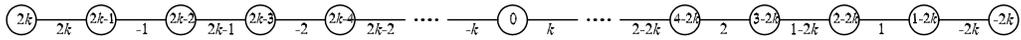}
\caption{A super edge-graceful labeling $P_{n}$ for $n\equiv 1 (\mod
4)$.}\label{1mod4}
\end{center}
\end{figure}

Notice that $P_{4k+1}$ consists of two paths $P_{2k+1}$
and $P_{2k+1}'$ with  edge labelings $f$ and $f'$, respectively, such that $f$ is inverse of $f'$.\\

{\noindent \emph{Case 2.}} $n\equiv 3 (\mod 4)$.

It follows that $n=4k+3$ where $k$ is a positive integer. Similarly
as in the previous case we  find a super edge-graceful labeling.
Notice that $P_{4k+3}$ consists of two paths $P_{2k+2}$ and
$P_{2k+2}'$ with edge labelings $f$ and $f'$, respectively, such
that $f$ is inverse of $f'$.\\
\begin{figure}[h]
\begin{center}
\includegraphics[width=13.5cm]{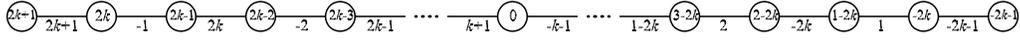}
\caption{A super edge-graceful labeling $P_{n}$ for $n\equiv 3 (\mod
4)$.}\label{3mod4}
\end{center}
\end{figure}

\begin{figure}[h]
\begin{center}
\includegraphics[width=7.6cm]{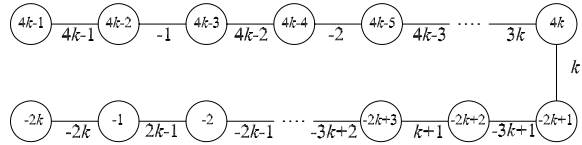}
\caption{A  labeling of $P_{4k}$.}\label{0mod8}
\end{center}
\end{figure}

 {\noindent \emph{Case 3.}} $n\equiv 0 (\mod 8)$.\\
 It follows that $n=8k$ for some positive integer $k$. We will consider $P_n$
 as a union of $P_{4k}$ and $P_{4k}'$ with  edge labelings $f$ and $f'$, respectively, such that $f$ is inverse of $f'$.
 We join the paths $P_{4k}$ and $P_{4k}'$ by an edge with label
 $0$. We label the edges along $P_{4k}$ by
 $(4k-1,-1,4k-2,-2,4k-3,-3,\ldots,-k+1,3k,k,-3k+1,k
 +1,-3k+2,k+2,\ldots,-2k-1,2k-1,-2k)$ (see Figure~\ref{0mod8}).

\begin{figure}[h]
\begin{center}
\includegraphics[width=7.6cm]{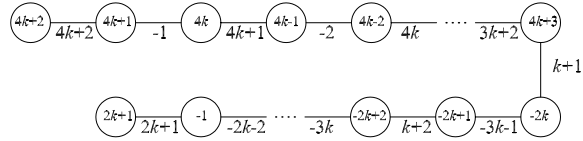}
\caption{A labeling of $P_{4k+3}$.}\label{6mod8}
\end{center}
\end{figure}

{\noindent \emph{Case 4.}} $n\equiv 6 (\mod 8)$.\\
Let  $n=8k+6$ for some positive integer $k$. As in Case 3  we will
consider $P_n$
 as a union of $P_{4k+3}$ and ${P}_{4k+3}'$ with  inverse labelings $f$ and $f'$, respectively, joined by an edge with label
 $0$. We label the edges along $P_{4k+3}$ by
 $(4k+2,-1,4k+1,-2,4k,-3,\ldots,3k+3,-k,3k+2,k+1,-3k-1,k
 +2,-3k,k+3,\ldots,2k,-2k-2,2k+1)$  (see Figure~\ref{6mod8}).

\begin{figure}[h]
\begin{center}
\includegraphics[width=9cm]{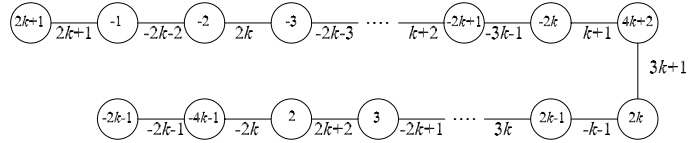}
\caption{A labeling of $Q_{1}$.}\label{4mod8_1}
\end{center}
\end{figure}

{\noindent \emph{Case 5.}} $n\equiv 4 (\mod 8)$.\\
Let  $n=8k+4$ for some positive integer $k$. We will consider $P_n$
 as a union of paths $Q_1$ and $Q_{2}$ of lengths $4k+3$ and $4k+1$, respectively, joined by an edge with label
 $0$.
 We label the edges along $Q_{1}$ by
 $(2k+1,-2k-2,2k,-2k-3,\ldots,k+2,-3k-1,k+1,3k+1,-k-1,3k,\ldots,-2k+1,2k+2,-2k,-2k-1)$  (see Figure~\ref{4mod8_1}).

\begin{figure}[h]
\begin{center}
\includegraphics[width=9cm]{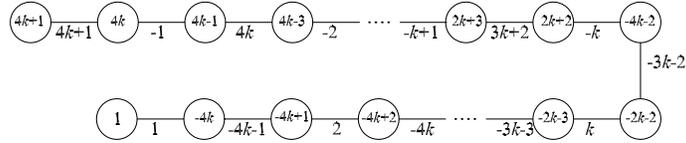}
\caption{A labeling of $Q_{2}$.}\label{4mod8_2}
\end{center}
\end{figure}

Further, we label edges along $Q_{2}$ by
 $(4k+1,-1,4k,-2,\ldots,-k+1,3k+2,-k,-3k-2,k,-3k-3,\ldots,-2k,2,-4k-1,1)$  (see Figure~\ref{4mod8_2}).
\\

\begin{figure}[h]
\begin{center}
\includegraphics[width=9cm]{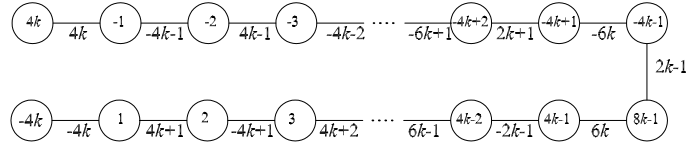}
\caption{A labeling of $Q_{3}$.}\label{2mod16_3}
\end{center}
\end{figure}

\begin{figure}[h]
\begin{center}
\includegraphics[width=13.7cm]{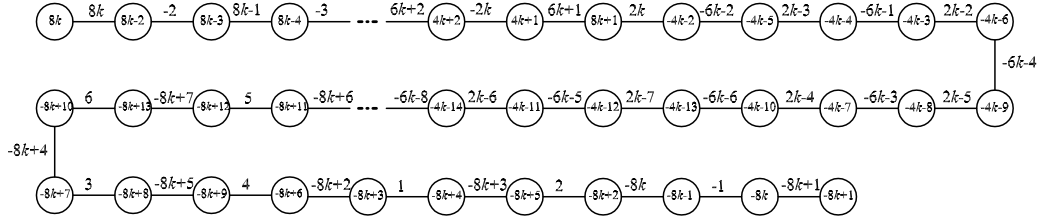}
\caption{A labeling of $Q_{4}$.}\label{2mod16_4}
\end{center}
\end{figure}

{\noindent \emph{Case 6.}} $n\equiv 2 (\mod 16)$.\\
Let  $n=16k+2$ for some positive integer $k$. As in Case 5 we
consider $P_n$
 as a union of paths $Q_3$ and $Q_{4}$ of lengths $8k+2$ and $8k$, respectively, joined by an edge with label
 $0$. We label the edges along $Q_{3}$ by
 $(4k,-4k-1,4k-1,-4k-2,\ldots,-6k+1,2k+1,-6k,2k-1,6k,-2k+1,6k-1,\ldots,4k+2,-4k+1,4k+1,-4k)$  (see Figure~\ref{2mod16_3}).
Then we label edges of  $Q_4$ by
$(8k,-2,8k-1,-3,\ldots,6k+2,-2k,6k+1,2k,-6k-2,2k-3,-6k-1,2k-2,-6k-4,2k-5,-6k-3,2k-4,-6k-6,2k-7,-6k-5,-2k+6,-6k-8,\ldots,
-8k+6,5,-8k+7,6,-8k+4,3,-8k+5,4,-8k+2,1,-8k+3,2,-8k,-1,-8k+1)$ (see
Figure~\ref{2mod16_4}).\\

\begin{figure}[h]
\begin{center}
\includegraphics[width=9cm]{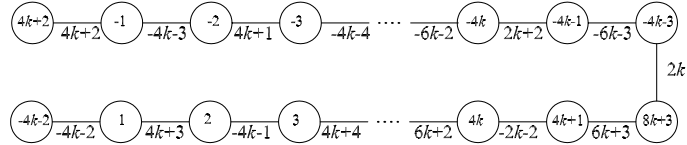}
\caption{A labeling of $Q_{5}$.}\label{10mod16_5}
\end{center}
\end{figure}

{\noindent \emph{Case 7.}} $n\equiv 10 (\mod 16)$.\\
 Let  $n=16k+10$ for some positive integer $k$. As in previous cases we  consider $P_n$
 as a union of paths $Q_5$ and $Q_{6}$ of lengths $8k+6$ and $8k+4$, respectively, joined by an edge with label
 $0$. We label the edges along $Q_{5}$ by
 $(4k+2,-4k-3,4k+1,-4k-4,-6k-2,2k+2,-6k-3,2k,6k+3,-2k+2,6k+2,\ldots,4k+4,-4k-1,4k+3,-4k-2)$  (see
Figure~\ref{10mod16_5}).\\

\begin{figure}[h]
\begin{center}
\includegraphics[width=13.7cm]{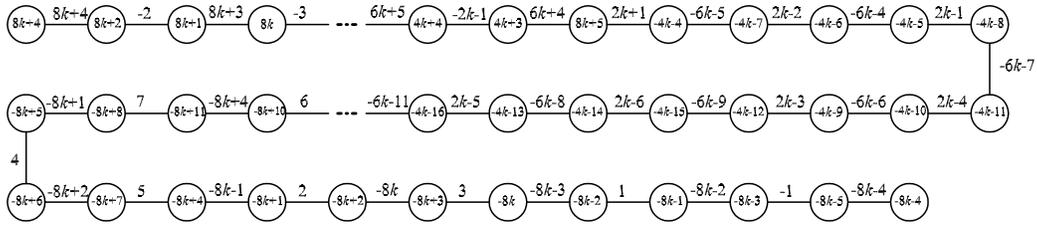}
\caption{A labeling of $Q_{6}$.}\label{10mod16_6}
\end{center}
\end{figure}

Then we label $Q_6$ by
$(8k+4,-2,8k+3,-3,\ldots,6k+5,-2k-1,6k+4,2k+1,-6k-5,2k-2,-6k-4,2k-1,-6k-7,2k-4,-6k-6,2k-3,-6k-9,2k-6,-6k-8,2k-5,-6k-11,\ldots,
6,-8k+4,7,-8k+1,4,-8k+2,5,-8k-12,-8k,3,-8k+3,1,-8k-2,-1,-8k-4)$
(see Figure~\ref{10mod16_6}). \qed\\

The corollary follows immediately from the proof
of~Theorem~\ref{CFX}.
\begin{cor}
If $n$ is odd then the cycle $C_n$ is super edge-graceful.
\end{cor}
\prf Since $n$ is odd then we can use the labeling for a path $P_n$,
and after that, by joining together the end vertices of the path by
the edge with label $0$, we obtain a cycle $C_n$ which is super
edge-graceful. \qed

\end{document}